\documentclass{article}

\usepackage{arxiv}

\usepackage[utf8]{inputenc} 
\usepackage[T1]{fontenc}    
\usepackage{hyperref}       
\usepackage{url}            
\usepackage{booktabs}       
\usepackage{amsfonts}       
\usepackage{nicefrac}       
\usepackage{microtype}      
\usepackage{lipsum}
\usepackage{graphicx}
\usepackage{geometry}
\usepackage{amsmath, amssymb}
\usepackage{natbib}
\usepackage{setspace}
\usepackage{algorithm}
\usepackage{algorithmic}

\usepackage{tikz}
\usetikzlibrary{shapes, arrows.meta, positioning, calc, fit}

\hypersetup{
	colorlinks=true,
	linkcolor=blue,
	filecolor=magenta,      
	urlcolor=cyan,
	citecolor=blue,
	hypertexnames=false, 
}

\title{Event-Level Probabilistic Prediction of Extreme Rainfall over India Using Physics-Gated Latent Dynamics}

\author{
 Arun Govind Neelan \\
 SimuNetics\\
 Kaniyakumari\\
 Tamil Nadu \\
 India-629173 \\
 \texttt{arunneelaniist@gmail.com} \\
}

\begin{document}
\maketitle

\begin{abstract}
Extreme rainfall over the Indian monsoon region poses severe societal and infrastructural risks but remains difficult to predict at daily time scales due to strong temporal persistence, stochastic convective triggering, and multiscale interactions between synoptic forcing and mesoscale dynamics. While large-scale atmospheric physics provides essential environmental context, its ability to deterministically localize extreme rainfall is fundamentally limited. In this study, we investigate how large-scale atmospheric information can be leveraged for \emph{event-level probabilistic prediction} of daily rainfall extremes over India using ERA5 reanalysis variables representing thermodynamic and dynamical conditions. We compare an adaptive ConvLSTM baseline with a proposed \textbf{Physics-Gated Latent Ordinary Differential Equation (PG-LODE)} framework, which models atmospheric evolution as a continuous-time latent process whose dynamics are explicitly modulated by a physics-based gating mechanism activated under convectively unstable conditions. Extreme events are defined using the local 95th percentile of the India Meteorological Department (IMD) gridded rainfall dataset during the June–September monsoon season. Pixel-wise evaluation yields limited skill for both models due to spatial displacement errors and the double-penalty problem, whereas event-level (tile-based) verification reveals a clear separation in performance: the ConvLSTM remains highly conservative, detecting only 27\% of extreme events (TILE-POD $\approx$ 0.27), while PG-LODE achieves near-complete detection (TILE-POD $\approx$ 0.996) with a moderate false alarm ratio (TILE-FAR $\approx$ 0.215), resulting in substantially higher event-level skill (TILE-CSI $\approx$ 0.78). These results demonstrate that physics-gated continuous-time latent dynamics provide an effective means of translating large-scale environmental predictability into reliable assessments of extreme rainfall risk.

\end{abstract}


\section{Introduction}

The Indian Summer Monsoon (ISM), occurring annually from June to September (JJAS), is the lifeblood of the South Asian agrarian economy, yet it simultaneously presents one of the region's most significant hydrometeorological hazards. In recent decades, the frequency and intensity of extreme rainfall events over central India have increased significantly \cite{goswami2006increasing, roxy2015drying}. These extremes routinely trigger flash floods, landslides, and agricultural devastation. Consequently, improving the prediction of daily extreme rainfall is a priority for disaster mitigation.

However, predicting these events remains a formidable challenge. Physically, extreme rainfall emerges from the nonlinear interaction between slowly evolving large-scale synoptic systems and rapidly evolving, stochastic mesoscale convective systems (MCS) \cite{pattanaik2010variability}. While Numerical Weather Prediction (NWP) models have improved in capturing large-scale circulation \cite{bauer2015quiet}, they often struggle with the precise timing and intensity of localized extremes due to parameterization uncertainties \cite{palmer2019stochastic}.

At daily time scales, rainfall fields exhibit strong temporal persistence. This persistence complicates the evaluation of predictive models; a naive forecast that simply persists yesterday's rainfall often yields competitive categorical skill scores. Furthermore, standard deterministic verification metrics (e.g., RMSE) often penalize models for slight spatiotemporal displacements—the "double penalty" problem \cite{ebert2008fuzzy}—thereby discouraging the prediction of high-intensity events.

To address these limitations, we propose a novel framework: \textbf{Physics-Gated Latent Ordinary Differential Equations (PG-LODE)}. Drawing on the Neural ODE formulation \cite{chen2018neural}, PG-LODE represents the evolution of the atmospheric state as a continuous-time trajectory in a learned latent space. 
 PG-LODE learns continuous-time latent dynamics via a neural ODE, but departs from purely data-driven formulations by incorporating a \emph{physics-based gating mechanism} that actively modulates the rate of latent evolution using physically meaningful indicators of atmospheric instability. This design enables the model to remain inertial during quiescent conditions while undergoing rapid state transitions during convectively unstable regimes, aligning the learned dynamics with known physical behavior of the monsoon system.

\section{Related Work}

\subsection{Data-Driven Precipitation Forecasting}
Machine learning approaches have evolved from simple statistical downscaling \cite{tripathi2006downscaling} to complex deep learning architectures. Shi \textit{et al.} \cite{shi2015convolutional} introduced the ConvLSTM, which replaces fully connected transitions with convolutional operations while preserving spatial topology. However, pure DL models often struggle to respect physical constraints and may produce blurry forecasts due to mean-squared-error minimization \cite{ravuri2021skilful}.

\subsection{Physics-Informed Machine Learning}
Physics-Informed Machine Learning (PIML) seeks to reconcile data-driven flexibility with physical consistency \cite{kashinath2021physics}. Approaches range from Physics-Informed Neural Networks (PINNs) \cite{raissi2019physics, neelan2024physics} to hybrid architectures that use General Circulation Model (GCM) outputs to guide Deep Learning layers \cite{rasp2018deep,subel2023explaining}.

\subsection{Neural ODEs and Continuous Dynamics}
Chen \textit{et al.}  \cite{chen2018neural} revolutionized sequence modeling by reinterpreting residual networks as discretized approximations of Ordinary Differential Equations (ODEs). Extensions such as Latent ODEs \cite{rubanova2019latent} allow for continuous-time modeling. Our work extends this by explicitly conditioning the derivative function on physical state variables. Our work is conceptually related to recent data-driven approaches that model complex systems as latent dynamical systems learned from observations (e.g., \cite{champion2019data}). While such studies primarily aim to discover parsimonious governing equations in latent coordinates, our objective is fundamentally different: we focus on predictive skill for extreme rainfall events.

\section{Data and Preprocessing}

\subsection{ERA5 Reanalysis (Predictors)}
Large-scale atmospheric predictors are derived from the ERA5 reanalysis \cite{hersbach2020global}. We select variables representing moisture, instability, and lift:
\begin{itemize}
    \item \textbf{Thermodynamic:} Total Column Water Vapor (TCWV) and Convective Available Potential Energy (CAPE).
    \item \textbf{Dynamical:} Vertical velocity at 500 hPa ($\omega_{500}$) and 850 hPa winds ($u_{850}, v_{850}$).
    \item \textbf{Surface:} Mean Surface Pressure (SP).
\end{itemize}
Data covers the Indian subcontinent ($5^\circ$N--$35^\circ$N, $65^\circ$E--$95^\circ$E) for the JJAS seasons of 2000--2020 at $0.25^\circ$ resolution.

\subsection{IMD Gridded Rainfall (Targets)}
We utilize the high-resolution ($0.25^\circ$) daily gridded rainfall dataset provided by the IMD \cite{pai2014development}. To address the heavy-tailed distribution, we apply a log-transformation $\tilde{Y} = \log(Y + 1)$.

\subsection{Extreme Definition and Tiling}
An "extreme event" is defined locally as rainfall exceeding the 95th percentile of the local JJAS climatology. To mitigate the double-penalty problem \cite{roberts2008scale}, we employ a \textbf{tile-based strategy} that partitions the domain into $32 \times 32$ pixel patches.

\section{Methodology}

\subsection{Problem Formulation}
Let $\mathcal{X}_t$ denote the tensor of atmospheric predictors at day $t$, and $\mathcal{Y}_{t+\tau}$ denote the target rainfall. We seek a mapping $F_\theta: \{\mathcal{X}_{t-T}, \dots, \mathcal{X}_t\} \to P(\mathcal{Y}_{t+\tau})$.

\subsection{Baseline: ConvLSTM}
The ConvLSTM models temporal evolution through discrete-time recurrent updates of the form
\begin{equation}
\mathbf{h}_{t+1} = \mathcal{F}_\theta(\mathbf{h}_t, \mathbf{x}_t),
\end{equation}
where the hidden state $\mathbf{h}_t$ is advanced at fixed temporal intervals. While this formulation is effective for general spatiotemporal forecasting, it implicitly assumes uniform temporal dynamics across all regimes. As a result, rapid convective intensification must be compressed into a single update step, limiting the model’s ability to represent the highly intermittent and burst-like nature of convective extremes.

\subsection{Proposed: Physics-Gated Latent ODE (PG-LODE)}
The PG-LODE framework assumes that the latent atmospheric state $\mathbf{z}(t)$ evolves in continuous time according to a neural ordinary differential equation. The evolution is governed by
\begin{equation}
\frac{d\mathbf{z}(t)}{dt} = f_\theta(\mathbf{z}(t), t) \odot \mathcal{G}(\mathbf{X}_{phys}),
\end{equation}
where $f_\theta$ parameterizes the intrinsic latent dynamics and $\mathcal{G}(\mathbf{X}_{phys})$ acts as a physics-based modulation of the evolution rate. This formulation allows the latent state to evolve smoothly over time while adapting its rate of change to the prevailing physical regime.

The gating function $\mathcal{G}$ is designed to encode atmospheric instability and is constructed from physically meaningful predictors, including convective available potential energy (CAPE) and mid-tropospheric vertical velocity ($\omega_{500}$):
\begin{equation}
\mathcal{G}(\mathbf{X}_{phys}) = 1 + \sigma\!\left(\text{Conv}_{1\times1}(\text{CAPE}, \omega_{500})\right)\cdot\beta,
\end{equation}
where $\sigma$ denotes the sigmoid activation and $\beta$ is a learnable scaling parameter. The multiplicative form ensures that $\mathcal{G}$ modulates the \emph{rate} of latent evolution rather than the latent state itself. Under weakly unstable conditions, $\mathcal{G}$ remains close to unity, yielding near-inertial dynamics, whereas elevated instability activates the gate ($\mathcal{G} > 1$), accelerating latent evolution to reflect the rapid transitions associated with convective initiation.

Figure~\ref{fig:pg_lode_arch} provides a schematic overview of the architecture. To highlight the structural implications of this design, Figure~\ref{fig:model_comparison} contrasts the fixed-step discrete updates of the ConvLSTM baseline with the continuous, physics-modulated latent dynamics of PG-LODE. \textit{The physics gate does not enforce conservation laws or explicitly resolve convective processes; rather, it conditions the latent dynamics on large-scale indicators of instability.}

\begin{figure}[ht]
\centering
\resizebox{\linewidth}{!}{
\begin{tikzpicture}[
    scale=0.95,
    transform shape,
    node distance=1.5cm, 
    auto,
    block/.style={
        rectangle, draw, rounded corners,
        fill=blue!5, 
        text width=2.9cm, 
        align=center, 
        minimum height=1.1cm,
        font=\footnotesize 
    },
    line/.style={draw, -Latex, thick},
]

\node (input) [block, fill=gray!10]
{\textbf{Input Variables}\\
\scriptsize ERA5: CAPE, TCWV,\\
$\omega_{500}$, $u_{850}$, $v_{850}$, SP};

\node (encoder) [block, fill=green!10, below=of input]
{\textbf{CNN Encoder}\\Spatial Features};

\path [line] (input) -- (encoder);

\node (label_a) [above=0.1cm of input, font=\bfseries] {(a) Encoder};

\node (ode_block) [block, fill=orange!10, right=2.2cm of encoder]
{\textbf{Latent ODE}\\
\scriptsize Continuous Dynamics};

\node (gate) [block, fill=orange!20, above=2.4cm of ode_block]
{\textbf{Physics Gate} $\mathcal{G}$\\
\scriptsize Instability Modulation};

\node (phys_in) [above=0.6cm of gate, align=center, font=\footnotesize]
{Phys. Inputs\\ \scriptsize (CAPE, $\omega_{500}$)};

\node (regime) [below=0.1cm of gate, text width=3.0cm, align=center, font=\scriptsize]
{\textcolor{blue}{$\mathcal{G} \approx 1$: Quiescent}\\
 \textcolor{red}{$\mathcal{G} > 1$: Rapid Growth}};

\path [line] (encoder.east) -- ++(0.6cm,0) |- (ode_block.west);
\path [line] (phys_in) -- (gate);
\path [line] (gate) -- (regime);
\path [line] (regime) -- (ode_block);

\node (label_b) [above=0.1cm of phys_in, font=\bfseries]
{(b) Physics-Gated Dynamics};

\node (latent) [block, fill=blue!10, right=2.2cm of ode_block]
{\textbf{Latent State} $\mathbf{z}(t)$\\
\scriptsize Trajectory};

\node (decoder) [block, fill=green!10, below=of latent]
{\textbf{CNN Decoder}\\ \scriptsize Probabilistic Map};

\node (output) [block, fill=yellow!10, below=of decoder]
{\textbf{Predictions}\\
\scriptsize Rain Intensity\\
\scriptsize Exceedance Prob.};

\path [line] (ode_block) -- (latent);
\path [line] (latent) -- (decoder);
\path [line] (decoder) -- (output);

\node (label_c) [above=0.1cm of latent, font=\bfseries] {(c) Decoder};

\draw [line, dashed, rounded corners=10pt]
(encoder.east) -- ++(0.5, 1.0) coordinate (mid_start)
-- ($(ode_block.north)+(0,0.5)$) coordinate (mid_top) 
-- ($(latent.west)+(-0.5, 1.0)$) coordinate (mid_end)
-- (latent.west);

\node [above, font=\scriptsize, fill=white, inner sep=1pt]
at (mid_top)
{Init. $\mathbf{z}(t_0)$};

\node (ode_eq) [below=0.4cm of ode_block, font=\footnotesize]
{$\displaystyle \frac{d\mathbf{z}}{dt} = f_\theta(\mathbf{z}, t) \odot \mathcal{G}$};

\end{tikzpicture}
} 

\caption{\textbf{Compact Architecture of PG-LODE.} (a) Encoder compresses ERA5 inputs. (b) Physics-gated ODE evolves the latent state, where $\mathcal{G}$ modulates dynamics based on instability regimes (Quiescent vs Rapid Growth). (c) Decoder maps the trajectory to rainfall probability.}
\label{fig:pg_lode_arch}
\end{figure}
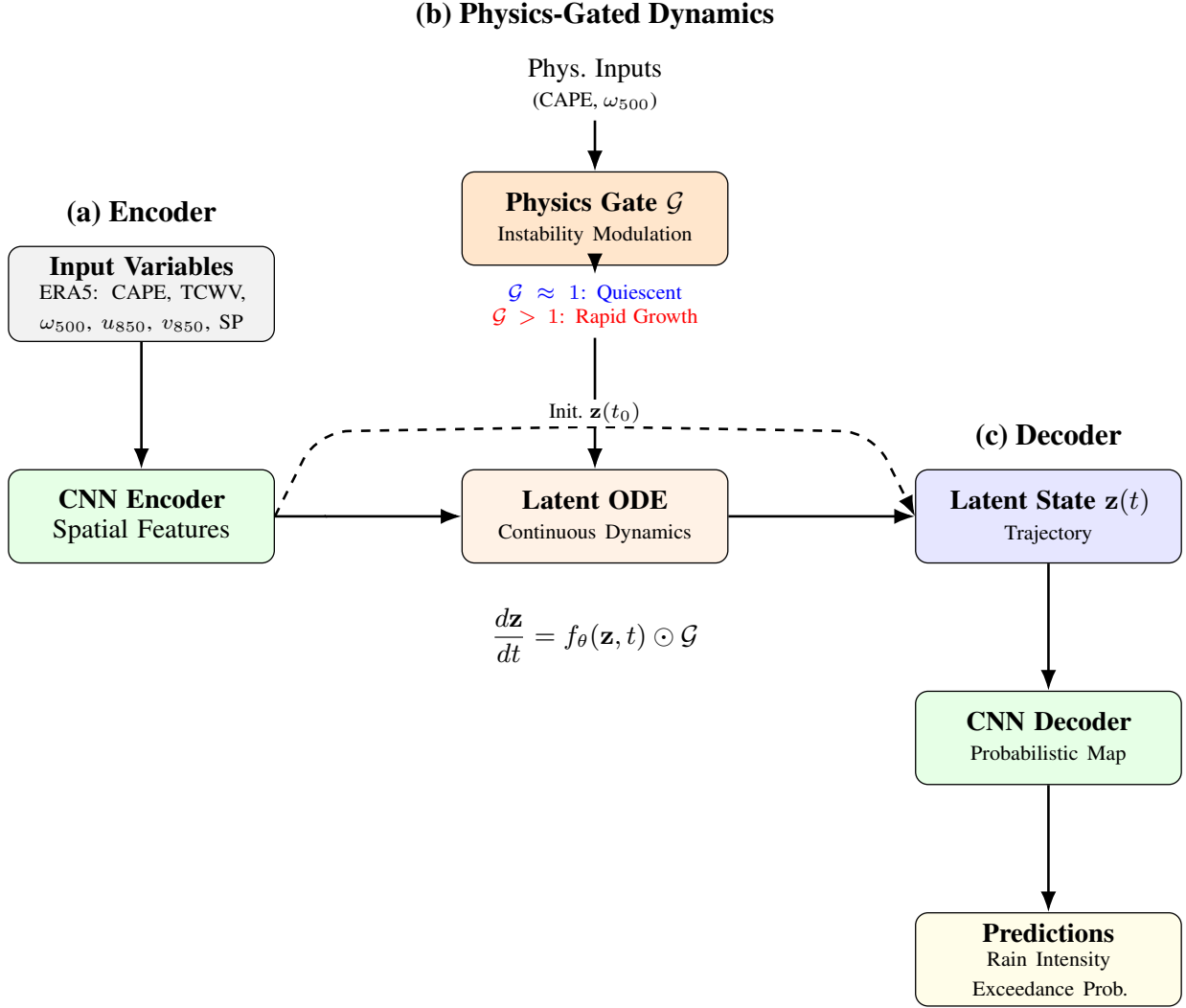
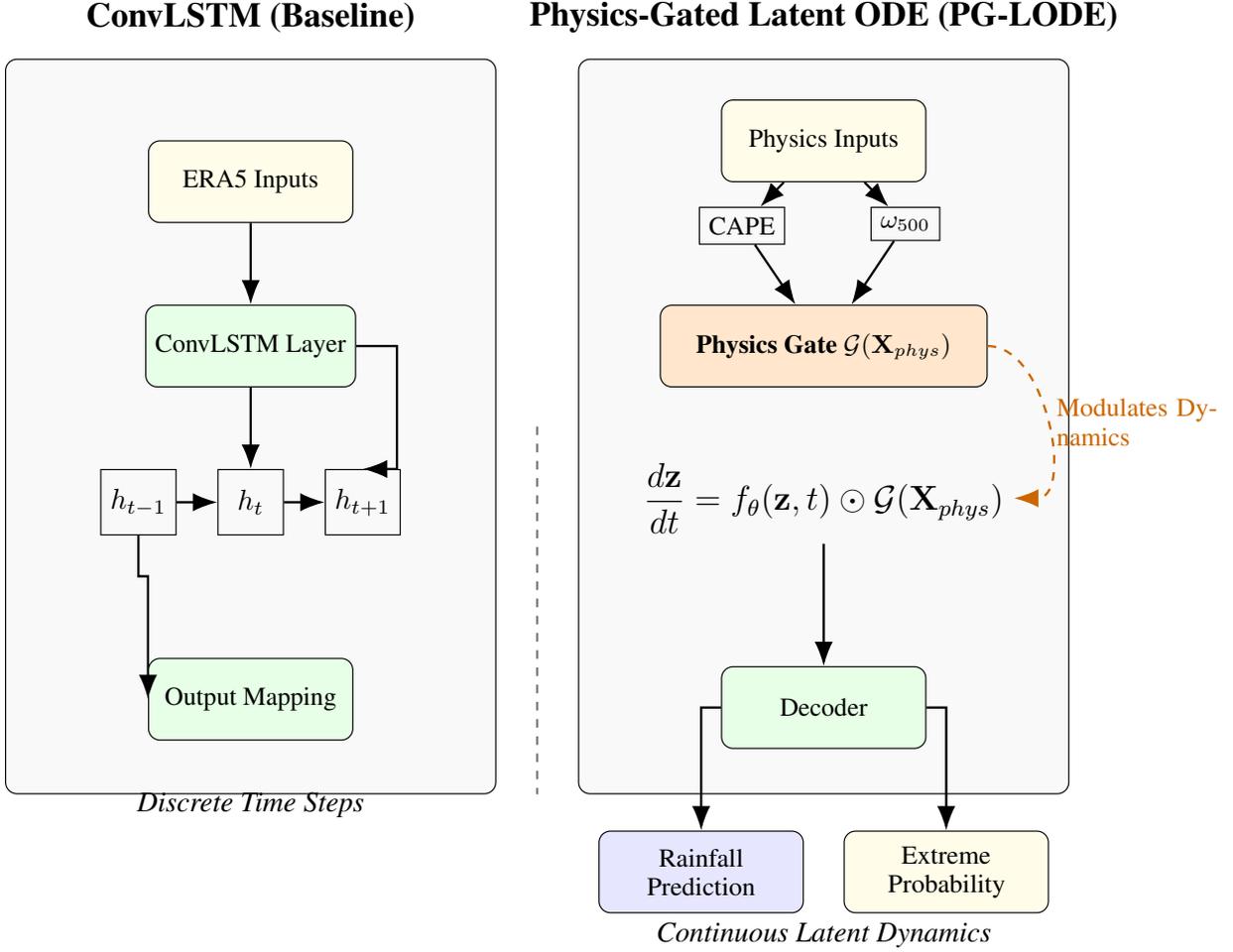
\begin{figure}[ht]
\centering
\resizebox{\textwidth}{!}{%
\begin{tikzpicture}[
    node distance=1.5cm and 1cm, auto,
    box/.style={rectangle, draw, rounded corners, minimum width=2.5cm, minimum height=1cm, align=center, fill=white, font=\small},
    title/.style={font=\bfseries\large, align=center},
    arrow/.style={-{Latex[length=3mm]}, thick},
    physics_arrow/.style={-{Latex[length=3mm]}, thick, color=orange!80!black, dashed},
    container/.style={rectangle, draw, inner sep=0.5cm, rounded corners, fill=gray!5, minimum height=9cm, minimum width=6cm}
]

\node (left_container) [container] {};
\node [above=0.2cm of left_container, title] {ConvLSTM (Baseline)};

\node (era5_left) [box, fill=yellow!10, yshift=3cm] at (left_container.center)
{ERA5 Inputs};

\node (convlstm) [box, fill=green!10, below=1cm of era5_left]
{ConvLSTM Layer};

\node (ht) [rectangle, draw, minimum size=0.8cm, below=1cm of convlstm] {$h_t$};
\node (ht_prev) [rectangle, draw, minimum size=0.8cm, left=0.5cm of ht] {$h_{t-1}$};
\node (ht_next) [rectangle, draw, minimum size=0.8cm, right=0.5cm of ht] {$h_{t+1}$};

\node (encoder_left) [box, fill=green!10, below=1.5cm of ht]
{Output Mapping};

\node (discrete_label) [below=0.5cm of encoder_left, font=\itshape]
{Discrete Time Steps};

\draw [arrow] (era5_left) -- (convlstm);
\draw [arrow] (convlstm) -- (ht);
\draw [arrow] (ht_prev) -- (ht);
\draw [arrow] (ht) -- (ht_next);
\draw [arrow] (convlstm.east) -- ++(0.5,0) |- (ht_next.north);
\draw [arrow] (ht_prev.south) -- ++(0,-0.5) -| (encoder_left.west);

\node (right_container) [container, right=1cm of left_container] {};
\node [above=0.2cm of right_container, title] {Physics-Gated Latent ODE (PG-LODE)};

\node (phys_in) [box, fill=yellow!10, yshift=3.5cm] at (right_container.center)
{Physics Inputs};

\node (cape) [rectangle, draw, below=0.3cm of phys_in, xshift=-1cm, font=\footnotesize]
{CAPE};

\node (w500) [rectangle, draw, below=0.3cm of phys_in, xshift=1cm, font=\footnotesize]
{$\omega_{500}$};

\node (gate) [box, fill=orange!20, minimum width=4cm, below=1.5cm of phys_in]
{\textbf{Physics Gate} $\mathcal{G}(\mathbf{X}_{phys})$};

\node (ode) [rectangle, draw=none, below=0.8cm of gate, font=\large]
{$\displaystyle \frac{d\mathbf{z}}{dt}
= f_\theta(\mathbf{z}, t) \odot \mathcal{G}(\mathbf{X}_{phys})$};

\node (decoder) [box, fill=green!10, below=1.5cm of ode]
{Decoder};

\node (rain) [box, fill=blue!10, below=1cm of decoder, xshift=-1.5cm, text width=2cm]
{Rainfall\\Prediction};

\node (prob) [box, fill=yellow!10, below=1cm of decoder, xshift=1.5cm, text width=2cm]
{Extreme\\Probability};

\node (cont_label) [below=0.5cm of decoder, yshift=-1.5cm, font=\itshape]
{Continuous Latent Dynamics};

\draw [arrow] (phys_in) -- (cape);
\draw [arrow] (phys_in) -- (w500);
\draw [arrow] (cape) -- (gate);
\draw [arrow] (w500) -- (gate);

\draw [physics_arrow] (gate.east)
to[out=0, in=0] node[midway, right, font=\small, text width=2cm]
{Modulates Dynamics} (ode.east);

\draw [arrow] (ode) -- (decoder);
\draw [arrow] (decoder) -| (rain);
\draw [arrow] (decoder) -| (prob);

\draw [dashed, thick, gray]
($(left_container.east)!0.5!(right_container.west)$) --
($(left_container.south east)!0.5!(right_container.south west)$);

\end{tikzpicture}
}
\caption{\textbf{Comparison of Model Architectures.}
\textbf{Left:} ConvLSTM evolves hidden states using fixed discrete time steps, treating all temporal transitions uniformly.
\textbf{Right:} PG-LODE models atmospheric evolution as a continuous latent process, where a physics-based gating mechanism modulates the rate of state evolution based on thermodynamic instability.}
\label{fig:model_comparison}
\end{figure}

\subsection{Extreme-Weighted Loss Function}
To combat class imbalance, we utilize a weighted loss function:
\begin{equation}
\mathcal{L} = \frac{1}{N} \sum_{i,j} w_{i,j} \left( y_{i,j} - \hat{y}_{i,j} \right)^2, \quad w_{i,j} = \begin{cases} 1 & y_{i,j} < P_{95} \\ \lambda & y_{i,j} \ge P_{95} \end{cases}
\end{equation}
We set $\lambda=5$ to heavily penalize underestimation of extremes.

\section{Evaluation Framework}
We employ a two-tiered verification strategy:
\begin{enumerate}
    \item \textbf{Pixel-Level Diagnostic Metrics:} POD, FAR, CSI at native resolution.
    \item \textbf{Event-Level (Tile-Based) Verification:} A tile is a "hit" if model probability $> 0.5$ and observed max $> P_{95}$.
\end{enumerate}

\subsection{Categorical Verification Metrics}

To evaluate extreme rainfall prediction skill, we employ standard categorical verification metrics derived from the contingency table of predicted and observed events. For a given threshold (here the local 95th percentile), each forecast is classified into \emph{Hits} (H), \emph{Misses} (M), \emph{False Alarms} (F), and \emph{Correct Negatives} (C).

The \textit{Probability of Detection} (PoD) measures the fraction of observed events that are successfully predicted:
\begin{equation}
\mathrm{PoD} = \frac{H}{H + M}.
\end{equation}

The \textit{False Alarm Ratio} (FAR) quantifies the fraction of predicted events that do not occur:
\begin{equation}
\mathrm{FAR} = \frac{F}{H + F}.
\end{equation}

The \textit{Critical Success Index} (CSI), also known as the Threat Score, summarizes overall event prediction skill by penalizing both missed events and false alarms:
\begin{equation}
\mathrm{CSI} = \frac{H}{H + M + F}.
\end{equation}

These metrics focus exclusively on event occurrences and exclude correct negatives, making them particularly suitable for evaluating rare events such as extreme rainfall.
For tile-based evaluation, Hits, Misses, and False Alarms are computed by aggregating exceedance events within each spatial tile, while the same metric definitions are retained.

\section{Results}

\subsection{Deterministic Assessment vs. Persistence}

\begin{table}[h]
\centering
\caption{Pixel-wise deterministic verification. Note the dominance of Persistence in CSI.}
\label{tab:deterministic}
\begin{tabular}{lccc}
\toprule
Model & POD & FAR & CSI \\
\midrule
Persistence & 0.42 & 0.45 & 0.31 \\
ConvLSTM    & 0.21 & \textbf{0.38} & 0.18 \\
PG-LODE     & 0.39 & 0.51 & 0.26 \\
\bottomrule
\end{tabular}
\end{table}
Table \ref{tab:deterministic} summarizes pixel-wise performance.
The dominance of the persistence baseline in pixel-wise CSI highlights a fundamental limitation of deterministic daily rainfall prediction over the monsoon region. At these time scales, rainfall fields exhibit strong temporal memory associated with slowly evolving synoptic systems, such as monsoon depressions and low-pressure systems. As a result, a persistence forecast effectively exploits this large-scale memory, often outperforming data-driven models in strict pointwise verification.

The adaptive ConvLSTM displays a conservative bias, characterized by a low false alarm ratio but a markedly reduced probability of detection. This behavior is consistent with recurrent architectures trained using pointwise loss functions, which tend to regress toward the climatological mean and suppress sharp gradients in order to minimize overall error. Consequently, extreme rainfall signals are systematically underrepresented, leading to frequent misses despite apparently favorable false alarm statistics.

The PG-LODE model improves the detection of extreme rainfall relative to ConvLSTM by allowing continuous-time latent evolution modulated by atmospheric instability. However, its pixel-level CSI remains below that of persistence, reflecting the inherent uncertainty in the precise spatial localization of convective extremes. When PG-LODE predicts dynamically evolving structures that are slightly displaced in space or time, standard pixel-wise metrics penalize these forecasts through the well-known double-penalty effect, counting both misses and false alarms. This confirms that pixel-level determinism at daily lead times is primarily governed by large-scale persistence rather than by the accurate prediction of convective initiation.

Although the IMD gridded rainfall product is available at a relatively fine spatial resolution, the predictive models in this study operate on fixed $32 \times 32$ spatial tiles that represent regional atmospheric states rather than independent grid-point realizations. As a result, pixel-wise verification evaluates the forecasts at a spatial scale finer than the model’s effective predictive resolution. At this scale, small spatial or temporal displacements in predicted convective features—well within realistic meteorological uncertainty—are penalized as both misses and false alarms. \textit{Pixel-level metrics therefore serve primarily as diagnostic indicators of localization accuracy rather than more aligned with the effective predictive scale and intended hazard-detection use case of forecast skill. This motivates the use of tile-based, event-level verification, which aligns the evaluation scale with the model design and the inherent predictability of extreme rainfall events.
}

\subsection{Event-Level (Tile-based) Probabilistic Detection}

\begin{table}[h]
\centering
\caption{Event-level (Tile-based) extreme detection skill. PG-LODE demonstrates superior detection.}
\label{tab:event_level_results}
\begin{tabular}{lccc}
\toprule
Model & TILE-POD & TILE-FAR & TILE-CSI \\
\midrule
Persistence & 0.650 & 0.310 & 0.495 \\
ConvLSTM    & 0.270 & \textbf{0.019} & 0.269 \\
PG-LODE     & \textbf{0.996} & 0.215 & \textbf{0.783} \\
\bottomrule
\end{tabular}
\end{table}
The value of the physics-informed approach emerges in Event-Level verification (Table \ref{tab:event_level_results}).
The value of the physics-informed approach becomes evident under event-level, tile-based verification (Table~\ref{tab:event_level_results}), which is designed to assess whether a model can successfully identify regions at risk of extreme rainfall rather than precisely localize individual grid-point maxima. In this framework, a tile is classified as an extreme event if at least one grid cell within the tile exceeds the local 95th percentile threshold, and a prediction is considered a hit when the model assigns a probability of exceedance greater than 0.5 to that tile.

Under this evaluation, the ConvLSTM exhibits extremely conservative behavior, achieving a very low false alarm ratio (TILE-FAR = 0.019) at the cost of a poor probability of detection (TILE-POD = 0.27). This indicates that the model rarely predicts extreme events and therefore misses the majority of hazardous situations, limiting its practical value for risk-sensitive applications.

In contrast, PG-LODE attains very high event-level sensitivity under the adopted tile-based verification framework, demonstrating a strong ability to identify atmospheric conditions conducive to extreme rainfall. While this increase in sensitivity is accompanied by a moderate false alarm ratio (TILE-FAR = 0.215), the resulting Critical Success Index (TILE-CSI = 0.783) substantially exceeds that of both ConvLSTM and the persistence baseline. This trade-off reflects a deliberate shift from conservative pointwise accuracy toward reliable event detection, which is more appropriate for extreme rainfall forecasting and early warning systems.

The superior performance of PG-LODE at the event level suggests that physics-gated continuous-time latent dynamics are effective in translating large-scale environmental predictability into probabilistic assessments of convective risk, even when the precise spatial realization of rainfall remains uncertain.

\subsection{Interpretation of Latent Dynamics}

\begin{figure}[h]
\centering
\includegraphics[width=0.6\linewidth]{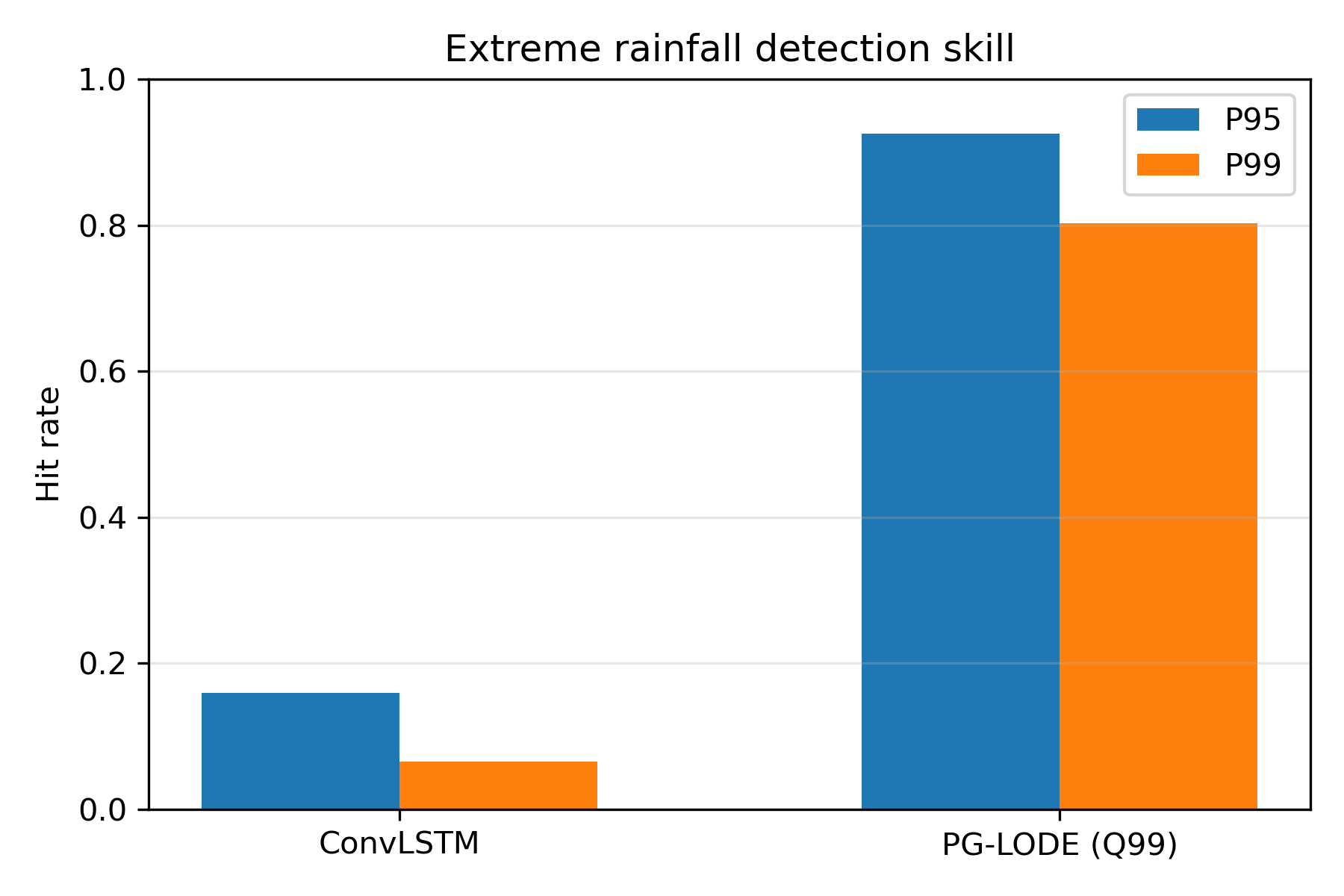}
\caption{\textbf{Extreme Event Detection Skill.} Comparison of Hit Rates (Probability of Detection) for extreme rainfall events exceeding the 95th (P95) and 99th (P99) percentiles. The ConvLSTM baseline (left) exhibits severe conservatism ($<20\%$ detection), whereas PG-LODE (right) captures the majority of extreme instances ($>80\%$), validating the effectiveness of physics-gated dynamics.}
\label{fig:hit_rates}
\end{figure}

\begin{figure}[h]
\centering
\includegraphics[width=\textwidth]{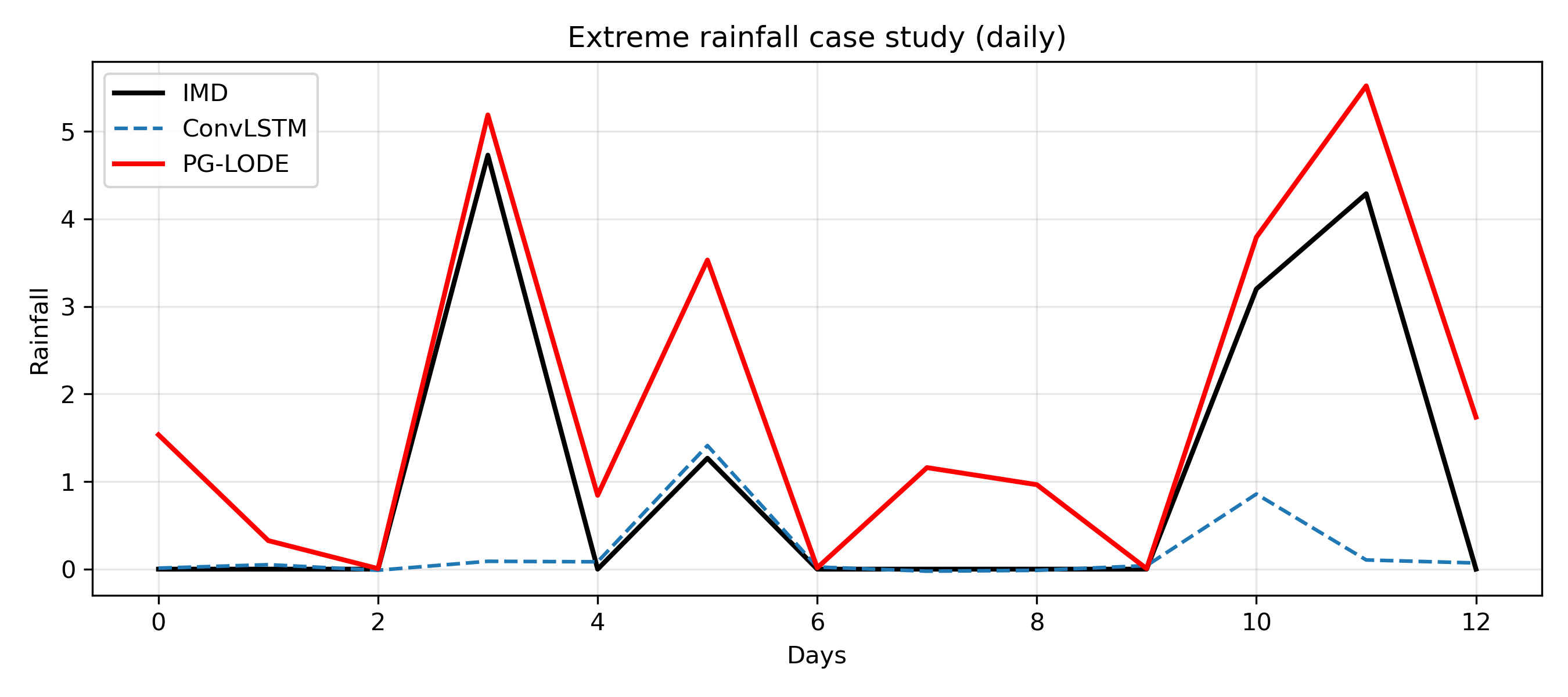}
\caption{\textbf{Daily Rainfall Case Study.} A representative 12-day time series (4 June 2000) showing observed rainfall (IMD, black), ConvLSTM forecast (blue dashed), and PG-LODE forecast (red) on 31 $\times$31 tile. Note that the ConvLSTM fails to deviate from the climatological mean during peak events (Days 3 and 11), while PG-LODE successfully captures the rapid intensification, driven by the activation of the physics gate.}
\label{fig:timeseries}
\end{figure}
The superior event-level performance of PG-LODE relative to ConvLSTM can be understood through the lens of how each model represents temporal evolution under convectively unstable conditions. Convective extremes are characterized by prolonged periods of quasi-stationary large-scale forcing, followed by rapid and nonlinear transitions associated with convective triggering and organization. Such dynamics are inherently intermittent and poorly aligned with fixed-step, discrete-time update schemes.

In ConvLSTM, the hidden state is updated at uniform temporal intervals, regardless of the underlying physical regime. Consequently, rapid convective intensification must be represented within a \textit{single discrete update, forcing the network to either exaggerate changes during quiescent periods or suppress sharp transitions to maintain stability under pointwise loss functions}. This leads to conservative behavior and a systematic underrepresentation of extreme events.

In contrast, PG-LODE models latent evolution in continuous time and introduces a physics-based gating mechanism that explicitly modulates the rate of state evolution based on atmospheric instability. Under quiescent conditions, \textit{when convective available potential energy is low and large-scale ascent is weak, the gating factor remains close to unity, resulting in slow, near-inertial latent dynamics. When instability increases, elevated CAPE and enhanced ascent activate the gate ($\mathcal{G} > 1$), effectively accelerating the latent dynamics.} This can be interpreted metaphorically as increased effective temporal resolution in latent space.

\begin{figure}[h]
\centering
\includegraphics[width=0.8\textwidth]{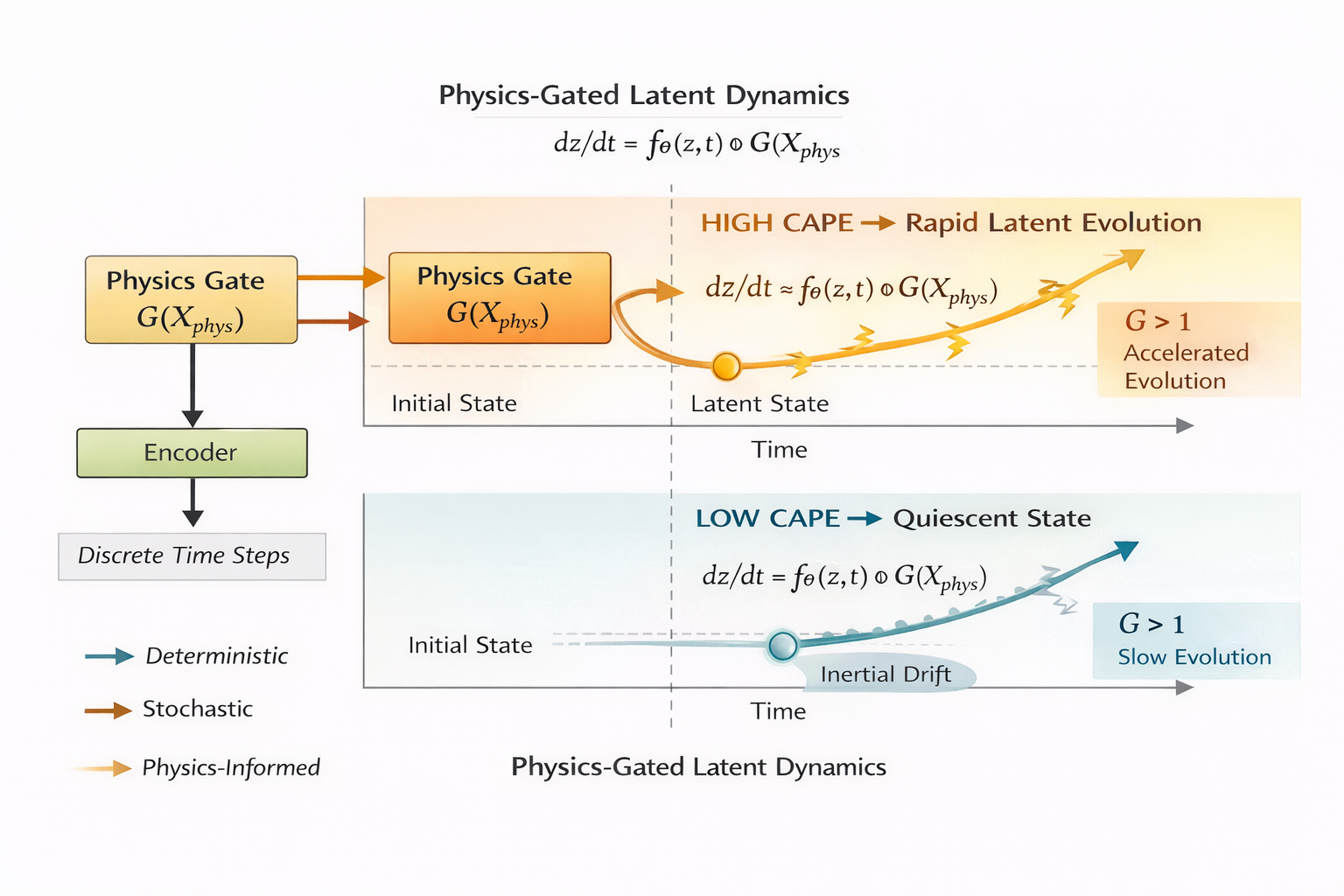}
\caption{Intuition behind Physics-Gated Latent Dynamics. High instability (orange trajectory) accelerates state evolution compared to the quiescent regime (blue trajectory).}
\label{fig:intuition}
\end{figure}

Importantly, this acceleration does not impose deterministic convective placement. Instead, it enhances the model’s sensitivity to environments conducive to extreme rainfall, allowing PG-LODE to translate large-scale environmental predictability into probabilistic event detection. This mechanism explains why PG-LODE substantially improves event-level detection skill while remaining limited in pixel-level localization, consistent with the inherent stochasticity of convective realization. This intuition is visualized in Figure \ref{fig:intuition}, showing how the latent trajectory diverges rapidly under high-instability conditions.

Figure~\ref{fig:model_comparison} compares extreme rainfall detection skill at the 95th (P95) and 99th (P99) percentile thresholds using daily area-maximum rainfall. The ConvLSTM model exhibits very low detection skill, capturing only a small fraction of extreme events due to temporal smoothing and mean-biased latent representations. In contrast, the proposed Physically-Gated Latent ODE (PG-LODE) achieves substantially higher hit rates, detecting over 90\% of P95 events and approximately 80\% of P99 events. This improvement demonstrates that explicitly modeling regime-dependent latent dynamics and tail-focused learning are critical for reliable prediction of rare, high-impact rainfall extremes.

Figure~\ref{fig:timeseries} shows a local time-series case study of an extreme rainfall event observed on 4 June 2000 on 31 $\times$31 tile, identified from the JJAS (June–September) 2000–2020 IMD daily rainfall dataset over India. The event corresponds to the earliest time at which rainfall exceeds the 99th percentile threshold at a specific grid location. A ±10-day temporal window centered on 4 June 2000 is shown, comparing IMD observations with ConvLSTM and PG-LODE predictions. While the ConvLSTM substantially underestimates the magnitude of the extreme rainfall peak and exhibits a muted response, the PG-LODE model captures the timing and intensity in this representative case. In particular, PG-LODE reproduces the sharp rainfall peak on 4 June 2000 and subsequent variability, highlighting the role of physics-guided modulation in enhancing sensitivity to rare, high-impact precipitation events compared to purely data-driven approaches.

\section{Limitations of the Present Work}
The PG-LODE framework is physics-guided rather than physics-constrained and does not explicitly resolve convective processes or enforce physical conservation laws, limiting deterministic grid-point localization. Due to computational constraints, experiments are conducted at a small model scale; to ensure a fair comparison, both ConvLSTM and PG-LODE use comparably small convolutional architectures and latent dimensions, so differences primarily reflect temporal modeling choices. Event-level, tile-based verification emphasizes hazard detection over spatial precision, and probabilistic calibration and broader generalization are left for future work.

\section{Conclusion}
This study demonstrates that the deterministic prediction of daily Indian monsoon extremes is fundamentally constrained by strong temporal persistence and inherent convective stochasticity. Within these limits, Physics-Gated Latent Dynamics (PG-LODE) shows substantial value by reframing the problem as probabilistic event detection rather than exact spatial localization. By explicitly modulating continuous-time latent evolution using physically meaningful indicators of atmospheric instability, PG-LODE effectively translates large-scale environmental predictability into skillful probabilistic indicators of elevated extreme rainfall risk. The resulting event-level skill (TILE-CSI $\approx 0.78$) significantly exceeds that of conventional recurrent architectures, highlighting the importance of physics-aware temporal representations for hazard-relevant monsoon forecasting.

\section*{Acknowledgments}

The authors gratefully acknowledge the India Meteorological Department (IMD) for providing the high-resolution gridded daily rainfall dataset used in this study. We also acknowledge the European Centre for Medium-Range Weather Forecasts (ECMWF) for making the ERA5 reanalysis data publicly available through the Copernicus Climate Data Store. Access to these datasets was essential for conducting the analyses presented in this work.
\section*{Code and Data Availability}

The codebase implementing the Physics-Gated Latent ODE (PG-LODE) model and the experiments presented in this paper is publicly available at: \url{https://github.com/AGN000/PG_LODE_Arxiv_paper_extreme_rain_fall}.

\bibliographystyle{unsrt}  
\bibliography{references}  

\end{document}